\documentclass[10pt,a4paper]{article}
\makeatletter
\usepackage[left=3cm,top=3cm,bottom=3cm,right=3cm,head=1cm,foot=1cm]{geometry}
\usepackage{amssymb}
\usepackage{amsmath}
\usepackage{amsthm}
\usepackage{amsfonts}
\usepackage{amssymb}
\usepackage{amscd}
\usepackage{url}
\usepackage[all]{xy}
\newcommand{\h}{\mathbf{H}}

\newcommand{\g}{\mathfrak{g}}

\newcommand{\5}{\hspace{0,5cm}}

\newcommand{\ad}{\stackrel{\mathrm{ad}}{\triangleright}}

\newcommand{\bad}{\stackrel{\mathrm{ad}}{\blacktriangleright}}
\newcommand{\obad}{\stackrel{\mathrm{oad}}{\blacktriangleright}}

\title{Covariant Dirac Operators on Quantum Groups}
\author{Antti J. Harju\footnote{ Department of Mathematics and Statistics, University of Helsinki, antti.harju@helsinki.fi}}
\date{}
\begin{document}
\maketitle

\begin{abstract} 
We give a construction of a Dirac operator on a quantum group based on any simple Lie algebra of classical type. The Dirac operator is an element in the vector space $U_q(\g) \otimes \mathrm{cl}_q(\g)$ where the second tensor factor is a $q$-deformation of the classical Clifford algebra. The tensor space $ U_q(\g) \otimes \mathrm{cl}_q(\g)$ is given a structure of the adjoint module of the quantum group and the Dirac operator is invariant under this action. The purpose of this approach is to construct equivariant Fredholm modules and $K$-homology cycles. This work generalizes the operator introduced by Bibikov and Kulish in \cite{BK}. \\


\end{abstract}

\textbf{1.} The Dirac operator on a simple Lie group $G$ is an element in the noncommutative Weyl algebra $U(\g) \otimes \mathrm{cl}(\g)$ where $U(\g)$ is the enveloping algebra for the Lie algebra of $G$. The vector space $\g$ generates a Clifford algebra $\mathrm{cl}(\g)$ whose structure is determined by the Killing form of $\g$. Since $\g$ acts on itself by the adjoint action, $ U(\g) \otimes \mathrm{cl}(\g)$ is a $\g$-module. The Dirac operator on $G$ spans a one dimensional invariant  submodule of $ U(\g) \otimes \mathrm{cl}(\g)$ which is of the first order in $\mathrm{cl}(\g)$ and in $U(\g)$. Kostant's Dirac operator \cite{Kos} has an additional cubical term in $\mathrm{cl}(\g)$ which is constant in $U(\g)$. 

The noncommutative Weyl algebra is given an action on a Hilbert space which is also a $\g$-module. Since the Dirac operator is the invariant subspace it follows that it commutes with the action of $\g$ and hence the Dirac operator acts as a constant on each irreducible component in the representation of $\g$ on the Hilbert space. The spectrum of the Dirac operator captures the metric properties of the the Riemannian manifold $G$ \cite{Con}. 

In \cite{BK} Bibikov and Kulish considered the quantum group deformation of $\mathfrak{su(2)}$ and constructed a Dirac operator which is invariant under the adjoint action of the quantum group.  In this approach the spectrum of the operator grows exponentially as a function of the highest weight of the representation of the quantum group. The sign operators of this type of Dirac operators can be used to define equivariant Fredholm modules which can be applied in the study of $D$-branes, $K$-homology, index theory and cyclic cohomology. 

The purpose of this paper is to study further the approach of the reference \cite{BK}. We are looking for an operator in the vector space $U_q(\g) \otimes \mathrm{cl}_q(\g)$, where $\mathrm{cl}_q(\g)$ is a $q$-deformation of the Clifford algebra which transforms covariantly under the action of the quantum group. Furthermore, we postulate the following defining principles for the covariant Dirac operator $D$ 
\begin{enumerate}
	\item Considered as a submodule of $U_q(\g) \otimes \mathrm{cl}_q(\g)$, $D$ spans a one dimensional trivial module under the adjoint action of $U_q(\g)$.
	\item Considered as an operator on a Hilbert space, $D$ commutes with the representation of $U_q(\g)$.
\end{enumerate}
For the Lie algebras the property 2. is a consequence of 1. However, if we let a quantum group act on the tensor product $U_q(\g) \otimes \mathrm{cl}_q(\g)$ with its coproduct, we need to choose the module structures in $\mathrm{cl}_q(\g)$ and $U_q(\g)$ carefully to make an operator with property 1. verify 2. This is not a general fact and will be explained in 8.

We consider the adjoint representation of the quantum group and define a bilinear form in this module which is invariant under the action. The braiding operator $\check{R}$ commutes with the coproduct of $U_q(\g)$ and it can be considered as a $q$-analogue of the permutation of a tensor product. We let $\check{R}$ act on a tensor product of adjoint representations and use the spectral decomposition of this action to define the $q$-Clifford algebra. The eigenvectors of $\check{R}$ split into two parts which can be considered as $q$-deformations of symmetric and antisymmetric tensor products. We identify the '$q$-symmetric' tensors with their image in the bilinear form. The practical difficulty in this approach is that there are no general formulas for the spectral decompositions. However, the explicit form of the $\check{R}$-operator is well known and so one can solve the eigenvalue problem with some mathematical software for any chosen quantum group and apply our results. 

We give a constructive proof for the existence of the covariant Dirac operator with the properties 1. and 2. on any quantum group based on any complex simple Lie algebra of classical type. The deformation parameter $q$ is supposed to be strictly positive real number. We give an explicit construction in the case of $SU_q(2)$ and build a Fredholm module from the sign operator. \\

\textbf{2. Conventions.}\ Let $\g$ be a simple finite dimensional Lie algebra with a set of simple roots $\Delta = \{\alpha_i: 1 \leq i \leq n\}$ and Cartan matrix $a_{ij}$. Let $q\neq 1$ be a complex number. The quantum group $U_q(\g)$ is the unital associative algebra with generators $k_{i}, k_{i}^{-1}, e_{i}, f_{i}$ ($1 \leq i \leq n$) subject to \cite{Dr,Ji} 
\begin{eqnarray*}
& &[k_{i}, k_{j}] = 0,\5  k_i k_i^{-1} = 1 \5 k_{i} e_j k_{i}^{-1} = q_i^{a_{ij}/2} e_j,\5 k_{i} f_j k_{i}^{-1} = q_i^{-a_{ij}/2} f_j, \\
& &[e_i, f_j] = \delta_{ij}\frac{k_{i}^2 - k_{i}^{-2}}{q_i -q_i^{-1}}, \\
& &\sum_{s=0}^{1-{a_{ij}}}(-1)^s \left[ \begin{array}{c}
1- a_{ij}  \\
s  
\end{array} \right]_{q_i} e_i^{s} e_j e_i^{1- a_{ij} - s} = 0 = \sum_{s=0}^{1-{a_{ij}}}(-1)^s \left[ \begin{array}{c}
1- a_{ij}  \\
s  
\end{array} \right]_{q_i} f_i^{s}f_j f_i^{1- a_{ij} - s}\5 (i \neq j),
\end{eqnarray*}
where $q_i = q^{d_i}$, $d_i$'s being the coprime integers such that $d_i a_{ij}$ is a symmetric matrix and the $q$-binomial coefficients are defined by
\begin{eqnarray*}
& &(m)_{q_i} = (q_i-q_i^{-1})(q_i^2-q_i^{-2}) \cdots (q_i^m -q_i^{-m}),\5 [0]_{q_i} =  1 \\
& &\left[ \begin{array}{c}
m  \\
n  
\end{array} \right]_{q_i} = \frac{(m)_{q_i}}{(n)_{q_i}(m-n)_{q_i}}.
\end{eqnarray*}
In the limit $q \rightarrow 1$ the algebra $U_q(\g)$ reduces to $U(\g)$. $U_q(\g)$ is a Hopf algebra with a coproduct $\triangle: U_q(\g) \rightarrow U_q(\g) \otimes U_q(\g)$, an antipode $S: U_q(\g) \rightarrow U_q(\g)$ and a counit $\epsilon: U_q(\g) \rightarrow \mathbb{C}$. The vectors $k_i$ and $k_i^{-1}$ are grouplike so that $\epsilon(k_i)= \epsilon(k_i^{-1}) = 1$ whereas $\epsilon(e_i) = \epsilon(f_i) = 0$. The coproduct $\triangle$ is noncocommutative and there exists a universal $R$-matrix such that 
\begin{eqnarray*}
R \triangle(x) R^{-1} = \sigma \triangle(x),
\end{eqnarray*}
where $\sigma$ permutes the tensor product. The $R$-matrix is an infinite sum defined in some completion of the tensor product $U_q(\g)\otimes U_q(\g)$ but only a finite number of terms are nonzero in any finite dimensional representation. Here we always assume that the parameter $q$ is a strictly positive real number. In this case the braiding operator $\check{R} = \sigma R$ is a selfadjoint operator in any finite dimensional representation. In the limit $q \rightarrow 1$, $\check{R}$ becomes the permutation operator. 

Let $\mathfrak{h}$ be the Cartan subalgebra of $\g$ with a basis $h_i$ satisfying $\alpha_j(h_i) = a_{ij}$. The representation $(V_{\lambda}, \pi_{\lambda,q})$ of $U_q(\g)$ is of the highest weight $\lambda \in \mathfrak{h}^*$ if there exists a vector $\xi$ so that $\pi_{\lambda,q} (e_i) \xi = 0$ and $\pi_{\lambda,q}(k_i)\xi = q^{\lambda(h_i)/2}\xi$ for all $i$ and the action of the operators $\pi_{\lambda,q}(f_i)$ on $\xi$ generate $V_{\lambda}$.

Denote by $P_+$ the set of integral dominant weights of a simple Lie algebra $\g$. As was shown in \cite{Lus,ros} the theory of finite dimensional representations of $\g$ and quantum group $U_q(\g)$ are identical in the case $q$ is not a root of unity. A highest weight module $(V_{\lambda},\pi_{\lambda,q})$ of $U_q(\g)$ is finite dimensional if and only if $\lambda \in P_+$. If $\lambda \in P_+$ then the dimension of each weight space is equal to the dimension of the corresponding weight space in the highest weight module $( V_{\lambda}, \pi_{\lambda})$ of $U(\g)$. The category of representations is semisimple with simple objects $(V_{\lambda}, \pi_{\lambda,q})$, $\lambda \in P_+$. For each morphism in the category of $U(\g)$ representations there exists a corresponding morphism in the category of $U_q(\g)$ representations.  

The matrix elements of an irreducible module $(V_{\lambda},\pi_{\lambda,q})$ of $U_q(\g)$ depend continuously on the parameter $q$ so that in the limit $q \rightarrow 1$, $\pi_{\lambda,q}(e_i)$ and $\pi_{\lambda,q}(f_i)$ define representation matrices for generators of $U(\g)$. These operators fix the representation $(V_{\lambda},\pi_{\lambda})$ of $U(\g)$ completely.   

We shall drop the $q$-subscript from the representations most often used in this work and denote by $(U, \pi)$, $(U^*, \pi^*)$ and $(V, \rho)$ the defining representation, its dual and the adjoint representation of $U_q(\g)$. 

The $q$-integers are defined by
\begin{eqnarray*}
[n] = \frac{q^n - q^{-n}}{q-q^{-1}}.
\end{eqnarray*}

\textbf{3.} The adjoint action of the quantum group on itself is an algebra homomorphism $U_q(\mathfrak{g}) \times U_q(\mathfrak{g}) \rightarrow U_q(\mathfrak{g})$ defined by
\begin{eqnarray*}
x \bad y = x'y S(x''),
\end{eqnarray*}
where $\triangle(x) = x' \otimes x''$. One can find a finite dimensional submodule in $U_q(\mathfrak{g})$ which is isomorphic to the adjoint representation of the quantum group \cite{Del}. Let $q = e^h$, $X = h^{-1} (R^{t}R-1)$ where $R^t = \sigma R \sigma$ and 
\begin{eqnarray*}
X_{lk} = (\pi_{lk} \otimes \mathrm{id})X \in \mathbb{C} \otimes U_q(\mathfrak{g}) \simeq  U_q(\mathfrak{g}),
\end{eqnarray*}
where $\pi_{lk}$ are the matrix elements of the defining representation of $U_q(\mathfrak{g})$. Denote by $\pi^*$ the dual representation 
\begin{eqnarray*}
\pi_{il}^*(x) = \pi_{li}(S(x)).
\end{eqnarray*}
The vectors $X_{lk}$ transform covariantly under the adjoint action 
\begin{eqnarray*}
x \bad X_{lk} = \sum_{i,j} X_{ij} \pi^*_{il}(x') \pi_{jk}(x''),\5 \hbox{for all}\5 x \in U_q(\mathfrak{g}).
\end{eqnarray*}

According to the termionology of \cite{Del}, a (weak) quantum Lie algebra is an invariant submodule in $U_q(\mathfrak{g})$ which is a deformation of $\g$ and transforms covariantly under the adjoint action. Denote by $\{u_i\}$ and $\{u_i^*\}$ the basis vectors of the defining representation and its dual and by $\{v_i\}$ the basis of the adjoint representation. Using the matrix coefficients of the module isomorphism $v_a \mapsto \sum_{i,j} K_a^{ij}(u_i^* \otimes u_j)$ we define
\begin{eqnarray*}
X_a = \sum_{i,j} K_{a}^{ij}(\pi_{ij} \otimes \mathrm{id}) X. 
\end{eqnarray*}
$X_a$'s span a quantum Lie algebra $\mathfrak{L}_q(\g)$ inside $U_q(\g)$ which is isomorphic to the adjoint representation of $U_q(\g)$. \\

\textbf{4.} Let $\g$ be a simple Lie algebra of classical type.  The adjoint representation can be considered as an invariant submodule $V \subset U^* \otimes U$ for the action
\begin{eqnarray*}
x \ad (u_l^* \otimes u_k) = \sum_{i,j} \pi_{li}(S(x')) u_i^* \otimes \pi_{jk}(x'') u_j 
\end{eqnarray*}
for all $x \in U_q(\g)$ and $u_l^* \otimes u_k \in V$. \\

\noindent \textbf{Proposition.} There exists a nondegenerate bilinear form $B_q: V \otimes V \rightarrow \mathbb{C}$ which is invariant under the adjoint action of $U_q(\g)$, i.e. 
\begin{eqnarray*}
B_q(x \ad (v \otimes w)):= B_q(x'\ad v \otimes x''\ad w) = \epsilon(x) B_q(v \otimes w). 
\end{eqnarray*}
$B_q$ is unique up to a multiplicative constant. \\
 
\noindent Proof. The adjoint representation $(V, \rho)$ and the dual representation $(V^*, \rho^*)$ are isomorphic. Let $\tau: V \rightarrow V^*$ denote a module isomorphism. Choose a basis $\{v_i\}$ of $V$ and let $\{v_i^*\}$ denote the dual basis. The canonical pairing defined on the generators by $\mathrm{eval}(v_j^* \otimes v_k) = v_j^* (v_k) = \delta_{jk}$ is nondegenerate. Thus, the composition
\begin{eqnarray*}
B_q:& & V \otimes V \stackrel{\sim}{\rightarrow} V^* \otimes V \stackrel{\mathrm{eval}}{\rightarrow} \mathbb{C}, \\
 & &v \otimes w \mapsto \mathrm{eval}(\tau(v) \otimes w)
\end{eqnarray*}
is nondegenerate. $B_q$ is invariant because
\begin{eqnarray*}
& &B_q(x'\ad v \otimes x''\ad w) = \mathrm{eval}(\tau( x'\ad v) \otimes x''\ad w) = \tau(v)(S(x')x''\ad (w)) = \epsilon(x) B_q(v \otimes w),
\end{eqnarray*}
for all $x \in U_q(\g)$ and $v,w \in V$. 

Let $\phi$ be the map $V \rightarrow V^*$ which sends $v \in V$ to the functional $\phi(v) (w) = B_q(v \otimes w) \in V^*$. Using the Hopf algebra axioms we see that $\phi$ is a module homomorphism: 
\begin{eqnarray*}
\phi(x\ad v)(w) &=& \phi(x'\ad v)(\epsilon(x'')w) = \phi(x'\ad v) ((x''S(x'''))\ad w) \\
&=& B_q(x'\ad v \otimes x''\ad (S(x''')\ad w)) = \epsilon(x') B_q(v \otimes S(x'')\ad w) \\
&=& \phi(v)(S(x)\ad w)
\end{eqnarray*}
for all $x \in U_q(\g)$ and $v, w \in V$. Furthermore, $\phi$ is a module isomorphism because $B_q$ is nondegenerate. If $\theta$ is another module isomorphism $\theta : V \rightarrow V^*$ we can define a module isomorphism $\theta^{-1} \circ \phi: V \rightarrow V$ which commutes with the action of the quantum group. $\g$ is simple and so the adjoint module $V$ is irreducible and the uniqueness follows from Schur's lemma. $\square$ \\

In practical calculations the form $B_q$ is easiest to find by fixing the costants directly from the invariance condition. \\

\textbf{5.} Let $\check{R}_i=\sigma_i R_i$ ( $1 \leq i \leq N-1$) be the linear operator where $R_i$ is the $R$ matrix acting on the $i$'th and $(i+1)$'th component in the tensor product space $V^{\otimes N}$ and $\sigma_i$ permutes the tensor components. The braiding operator $\check{R}_i$ commutes with the action of $U_q(\g)$ on the tensor product and thus the eigenspaces of $\check{R}_i$ are invariant subspaces of $U_q(\g)$. $\check{R}_i$ is a selfadjoint operator and its eigenvalues are real. Furthermore, the eigenvalue of a nonzero eigenspace is not equal to zero for any $q > 0$ because $\check{R}_i$ is an isomorphism. Thus, the tensor product splits into parts consisting of the vectors with strictly positive eigenvalues and strictly negative eigenvalues for any allowed value of $q$. In the classical limit $q \rightarrow 1$ these eigenspaces become the symmetric and antisymmetric tensor products in the $i$'th and $(i+1)$'th component. 

Given a spectral resolution of $\check{R}_i$ denote by $\{ a_{i,k}: k \in I\}$ the negative eigenvalues and by $\{ b_{i,k}: k \in J\}$ the positive eigenvalues of $\check{R}_i$. The braiding operators form a generalized Hecke-algebra with relations
\begin{eqnarray*}
& &\check{R}_{i} \check{R}_{i+1} \check{R}_{i} = \check{R}_{i+1} \check{R}_{i} \check{R}_{i+1}\\
& &\check{R}_{i} \check{R}_{j} = \check{R}_{j} \check{R}_{i},  \\
& & \prod_{k \in I}(\check{R}_{i}- a_{i,k})\prod_{l \in J}(\check{R}_{i}-b_{i,l}) = 0
\end{eqnarray*}
for all $1 \leq i,j \leq N-1$ and $|i-j|>1$. 

Let $T(V)$ be the tensor algebra of $V$. We define the covariant Clifford algebra as a projection of $T(V)$ on the $q$-analogue of the antisymmetric tensor products by
\begin{eqnarray*}
\mathrm{cl}_q(\g) = T(V)/\mathfrak{I}
\end{eqnarray*}
where the ideal $\mathfrak{I}$ is defined by  
\begin{eqnarray*}
\mathfrak{I} = \{ (\mathrm{id} - B_q^{i})v : v \in \mathrm{Ker}(\check{R}_i - b_{i,k})\5 \hbox{for some}\5 i \in \mathbb{N}, k \in J \}.
\end{eqnarray*}
$B_q^{i}$ is the invariant pairing of $i$'th and $(i+1)$'th tensor component. $\mathrm{cl}_q(\g)$ transforms covariantly under the adjoint action of $U_q(\g)$ because $B^i_q$ is invariant and the operators $\check{R}_i - b_{i,k}$ commute with the action of $U_q(\g)$. Therefore, it is a $U_q(\g)$-module algebra. Let us denote by $\gamma_q: V \rightarrow \mathrm{cl}_q(\g)$ the canonical embeddings. \\

\textbf{6.} There exists a Lie algebra homomorphism $\widetilde{\mathrm{ad}}: \g \rightarrow \mathrm{cl}(\g)$ which satisfies the equivariance condition
\begin{eqnarray*}
\gamma([x,y]) = [\widetilde{\mathrm{ad}}(x),\gamma(y)],
\end{eqnarray*}
for all $x,y \in \g$ where $\gamma$ is the canonical embedding $\g \rightarrow \mathrm{cl}(\g)$. Once the irreducible representation $(\Sigma,s)$ of $\mathrm{cl}(\g)$ is given this fixes the representation of $\g$ on $\Sigma$. \\

\noindent \textbf{Proposition.} The algebras $\mathrm{cl}_q(\g)$ and $U_q(\g)$ have representations $s_q$ and $\sigma_q$ on $\Sigma$ which satisfy the $q$-deformed equivariance condition 
\begin{eqnarray}\label{comp}
s_q (\gamma_q(x \ad v)) = \sigma_q(x') s_q(\gamma_q(v)) \sigma_q(S(x''))
\end{eqnarray}
for any $x \in U_q(\g)$ and $v \in V$. The representation $s_q$ of $\mathrm{cl}_q(\g)$ is irreducible. \\  

\noindent Proof. Denote by $(\Sigma, s)$ an irreducible $\mathrm{cl}(\g)$ representation. The vector space $B(\Sigma)$ of all endomorphisms of $\Sigma$ is a $\g$-module equipped with the action given by commutator
\begin{eqnarray*}
x \triangleright T =  s( \widetilde{\mathrm{ad}}(x')) T s (\widetilde{\mathrm{ad}}(S(x''))).
\end{eqnarray*}
Let us denote by $\sigma_q$ the representation of $U_q(\g)$ on $\Sigma$ which corresponds to the representation $s(\widetilde{\mathrm{ad}})$ of $U(\g)$ in the classical limit. The action
\begin{eqnarray*}
x \triangleright_q T =  \sigma_q(x') T \sigma_q(S(x''))
\end{eqnarray*}
defines a representation of $U_q(\g)$ on $B(\Sigma)$. The module $B(\Sigma)$ is reducible and the decomposition to irreducible subspaces is the same as in the classical case. 

The proof is based on the following observation. If $X$ and $Y$ are $U(\g)$ or $U_q(\g)$ submodules of $B(\Sigma)$ then the space $XY$ spanned by the matrix products is also a $U(\g)$ or $U_q(\g)$ submodule of $B(\Sigma)$. This is true because by using the Hopf algebra properties we get 
\begin{eqnarray*}
x \triangleright MN = (x' \triangleright M)(x'' \triangleright N) \in XY 
\end{eqnarray*}
for all $x \in U(\g)$, $M \in X$ and $N \in Y$. The corresponding formula holds also for the quantum group action $\triangleright_q$. In both cases the matrix multiplication defines a module homomorphism $m: X \otimes Y \rightarrow XY$ where the action on the tensor product is given by composing the coproduct with $\triangleright$ or $\triangleright_q$.

By equivariance, the submodule $B(\Sigma)' = s(\gamma(\g))$ is isomorphic to the adjoint representation. The module $B(\Sigma)' \otimes B(\Sigma)'$ reduces to invariant symmetric and antisymmetric components. The multiplication restricted to the symmetric submodule is a module homomorphism getting values in the trivial module $\mathbb{C} \textbf{1}$
\begin{eqnarray*}
m: s(\widetilde{\mathrm{ad}}(x)) \otimes s(\widetilde{\mathrm{ad}}(y)) + s(\widetilde{\mathrm{ad}}(y)) \otimes s(\widetilde{\mathrm{ad}}(x)) \mapsto B(x,y)\textbf{1},
\end{eqnarray*}
for all $x,y \in \g$, where $B$ is proportional to the Killing form. 

Denote by $B(\Sigma)_q'$ the submodule of $B(\Sigma)$ which is isomorphic to the adjoint representation of $U_q(\g)$ and reduces to the representation $s(\gamma(\g))$ in the classical limit. The tensor product $B(\Sigma)_q' \otimes B(\Sigma)_q'$ decomposes again into invariant components which are symmetric and antisymmetric in the $q$-deformed sense. The multiplication restricted to the $q$-symmetric part is a module homomorphism and must get values in the trivial module $\mathbb{C} \textbf{1}$ because the representation theory corresponds to the classical one. Especially the homomorphism must be determined by the invariant bilinear form which is unique. Therefore we find a representation $s_q$ of $\mathrm{cl}_q(\g)$ on $\Sigma$ which satisfies \eqref{comp} with $\sigma_q$.

It remains to show irreducibility. If the dimension of $\g$ is $2l$ the irreducible representation $\mathrm{cl}(\g) \rightarrow B(\Sigma)$ is an algebra isomorphism. If the dimension is $2l+1$ there exists two nonisomorphic irreducible representations which give an algebra isomorphism $\mathrm{cl}(\g) \rightarrow B(\Sigma) \oplus B(\Sigma)$. Since $\gamma(\g)$ generates $\mathrm{cl}(\g)$ as an algebra, each endomorphism algebra $B(\Sigma)$ must be generated by the subspace $s(\gamma(\g)) = B(\Sigma)'$. In terms of representation theory this has the following interpretation:  Starting with $B(\Sigma)' \otimes B(\Sigma)'$ one can apply the multiplication homomorphism to construct more irreducible $U(\g)$ representations and after taking finite number of steps each irreducible component of $B(\Sigma)$ is found. Similarly, we can consider $B(\Sigma)$ as $U_q(\g)$-module and apply the multiplication process to $B(\Sigma)_q'$ which gives all the irreducible components of irreducible $U_q(\g)$ representations in $B(\Sigma)$, but then the module algebra generated by $B(\Sigma)'_q$ is the whole space $B(\Sigma)$ and therefore no nontrivial invariant subspaces exist. $\5 \square$ \\

\noindent \textbf{Corollary} The algebras $\mathrm{cl}_q(\g)$ and $\mathrm{cl}(\g)$ are isomorphic as associative algebras. \\

\noindent Proof. It follows from the above proof that the representations $\mathrm{cl}_q(\g) \rightarrow B(\Sigma)$ in the $2l$ dimensional case and $\mathrm{cl}_q(\g) \rightarrow B(\Sigma) \otimes B(\Sigma)$ in the $2l+1$ dimensioal case are onto. Since the space of symmetric tensor products preserve the dimensionality in the $q$-deformation we get $\mathrm{dim}(\mathrm{cl}_q(\g))=\mathrm{dim}(\mathrm{cl}(\g))$. Then $\mathrm{dim}(\mathrm{cl}_q(\g)) = \mathrm{dim}(B(\Sigma))$ or $\mathrm{dim}(\mathrm{cl}_q(\g)) = 2\mathrm{dim}(B(\Sigma))$ in $2l$ and $2l+1$ dimensions. Thus, the irreducible representations define algebra isomorphisms 
\begin{eqnarray*}
& &s_q: \mathrm{cl}_q(\g) \rightarrow B(\Sigma),\5 \hbox{dim$(\g) = 2l$} \\
& &s_q^1 \oplus s_q^2 : \mathrm{cl}_q(\g) \rightarrow B(\Sigma) \oplus B(\Sigma),\5 \hbox{dim$(\g) = 2l+1$} . 
\end{eqnarray*}
Especially, the algebras $\mathrm{cl}_q(\g)$ and $\mathrm{cl}(\g)$ are isomorphic. $\5 \square$ \\

\textbf{7.} We first define an invariant one dimensional subspace in the vector space $V \otimes V^*$ and then define $D$ in the image of this subspace in a module isomorphism from $V \otimes V^*$ to a submodule of $U_q(\g) \otimes \mathrm{cl}_q(\g)$. The module structure of $  U_q(\g) \otimes \mathrm{cl}_q(\g)$ is chosen so that $D$ commutes with the representation. \\

\noindent \textbf{Proposition.} Let $\{v_i \}$ denote the basis of $V$ and $\{v_i^*\}$ the dual basis. The vector $\Omega \in V \otimes V^*$ defined by
\begin{eqnarray*}
\Omega = \sum_{i} v_i \otimes v_i^*
\end{eqnarray*}
is invariant under the action of $U_q(\g)$. \\

\noindent Proof. For all $x \in U_q(\g)$ 
\begin{eqnarray*}
x \ad \Omega &=& \sum_{i,k,l} \rho_{ki}(x') v_k \otimes \rho^*_{li}(x'') v^*_l \\
&=& \sum_{i,k,l} \rho_{ki}(x')\rho_{il}(S(x''))v_k \otimes v^*_l  \\
&=& \sum_{k,l} \rho_{kl}(x' S(x''))v_k \otimes v^*_l \\
&=& \sum_{k,l} \epsilon(x) \delta_{kl} v_k \otimes v^*_l = \epsilon(x)\Omega.\5 \square
\end{eqnarray*}

\textbf{8.} Representations on tensor products are dependent on the choice of the Hopf algebra structure for the quantum group. If we fix a Hopf structure, then by definition, the adjoint action on a vector $Z \otimes \psi \in U_q(\mathfrak{\g}) \otimes \mathrm{cl}_q(\g)$ is
\begin{eqnarray}\label{naive}
x \ad (Z \otimes \psi) := x' \bad Z \otimes x'' \ad \psi.
\end{eqnarray}
Let $\phi': V \rightarrow \mathfrak{L}_q(\g)$ and $\tau: V^* \rightarrow V$ be module isomorphisms. We could try to define the Dirac operator by $A = (\phi' \otimes \gamma_q \circ \tau)(\Omega)$. This is certainly invariant under the action \eqref{naive} of the quantum group. Let us write $A = \sum_{i,j} \alpha_{ij} X_i \otimes \psi_j$ for some complex numbers $\alpha_{ij}$. 

Denote by $s_q$ the irreducible representation of $\mathrm{cl}_q(\g)$ and $\sigma_q$ the representation of $U_q(\g)$ on $\Sigma$ which satisfy \eqref{comp}. Let $U_q(\g)$ act on $V_{\lambda} \otimes \Sigma$ by $x \mapsto \pi_{\lambda,q}(x') \otimes \sigma_q(x'')$. Even though $A$ is invariant it fails to commute with the representation of $U_q(\g)$. This can be seen using $x = x'\epsilon(x'')$ and $\epsilon(x) = S(x')x''$ twice
\begin{eqnarray*}
xA.| \omega \rangle &=& (\pi_{\lambda,q}(x') \otimes \sigma_q(x''))(\sum_{i,j} \alpha_{ij} \pi_{\lambda,q}(X_i) \otimes s_q(\psi_j))| \omega \rangle \\ 
&=& (\sum_{i,j} \alpha_{ij} \pi_{\lambda,q}(x' \bad X_i) \otimes s_q(x''' \ad \psi_j))(\pi_{\lambda,q}(x'') \otimes \sigma_q(x''''))| \omega \rangle \\
&\neq& \epsilon(x')(\sum_{i,j} \alpha_{ij} \pi_{\lambda,q}(X_i) \otimes s_q(\psi_j))(\pi_{\lambda,q}(x'') \otimes \sigma_q(x'''))| \omega \rangle = Ax.| \omega \rangle
\end{eqnarray*}
for some $x\in U_q(\g)$ and $| \omega \rangle \in V_{\lambda} \otimes \Sigma$ because of the noncocommutativity of the coproduct. \\

\textbf{9.} The problem above can be cured by modifying the structure of the module $\mathfrak{L}_q(\g)$. Let us choose the primary Hopf algebra structure for $U_q(\g)$ by
\begin{eqnarray*}
& &\triangle(e_i) = e_i \otimes k_i + k_i^{-1} \otimes e_i,\5 \triangle(f_i) = f_i \otimes k_i + k_i^{-1} \otimes f_i,\5 \triangle(k_i) = k_i \otimes k_i \\
& & S(k_i) = k_i^{-1},\5 S(e_i) = - q e_i,\5 S(f_i) = - q^{-1} f_i.
\end{eqnarray*}
The Sweedler's notation $\triangle(x) = x' \otimes x''$ is always applied to this one. The opposite Hopf algebra structure is defined by
\begin{eqnarray*}
\triangle^{\mathrm{op}}(x) = x'' \otimes x',\5 S^{\mathrm{op}}(k_i) = k_i^{-1},\5 S^{\mathrm{op}}(e_i) = - q^{-1} e_i ,\5 S^{\mathrm{op}}(f_i) = - q f_i.\
\end{eqnarray*}
We apply the method of 3. to construct a quantum Lie algebra for the Hopf algebra $U_q(\g)$ with the opposite Hopf structure. This gives us a module $\mathfrak{L}_q^{\mathrm{op}}(\g)$ with a basis $\{Z_i\}$ transforming as 
\begin{eqnarray*}
x \obad Z_i = x''Z_i S^{\mathrm{op}}(x') = \sum_j \rho_{ji}(x) Z_j
\end{eqnarray*}
for each $x \in U_q(\g)$ and $\rho_{ji}$ are the matrix elements of the fixed adjoint representation $V$. Denote by $\phi$ the module isomorphism $V \rightarrow \mathfrak{L}^{\mathrm{op}}_q(\g)$. \\

\noindent \textbf{Theorem.} Consider the representation theory of $\mathrm{cl}_q(\g)$ and $U_q(\g)$ given by $(\Sigma,s_q)$ and $\sigma_q$ as in the proposition of 6. The operator 
\begin{eqnarray*}
D = (\phi \otimes \gamma_q \circ \tau)(\Omega) \in \mathfrak{L}^{\mathrm{op}}_q(\g) \otimes \mathrm{cl}_q(\g) \subset U_q(\g) \otimes \mathrm{cl}_q(\g) 
\end{eqnarray*}
is invariant under the action of $U_q(\g)$:  $\triangle(x)( \obad \otimes \ad)D = \epsilon(x)D$. $D$ commutes with the action of $U_q(\g)$ on $V_{\lambda} \otimes \Sigma$ for any $\lambda \in P_+$.  \\

\noindent Proof. $D$ clearly spans a singlet in $U_q(\g) \otimes \mathrm{cl}_q(\g)$ because $\Omega$ does. Let us write $D = \sum_{i,j} \alpha_{ij} Z_i \otimes \psi_j$. Let $| \omega \rangle \in V_{\lambda} \otimes \Sigma$ and $x\in U_q(\g)$. Using the Hopf algebra properties
\begin{eqnarray*}
\epsilon(x) = S(x')x'' = S^{\mathrm{op}}(x'')x'
\end{eqnarray*}
we find that 
\begin{eqnarray*}
xD.| \omega \rangle &=& (\pi_{\lambda,q}(x') \otimes \sigma_q(x''))(\sum_{i,j} \alpha_{ij} \pi_{\lambda,q}(Z_i) \otimes s_q(\psi_j))| \omega \rangle \\
 &=& (\sum_{i,j} \alpha_{ij} \pi_{\lambda,q}(\epsilon(x')x'' Z_i)) \otimes (\sigma_q(x''' \epsilon(x'''') )  s_q(Z_j))| \omega \rangle \\
&=& (\sum_{i,j} \alpha_{ij} \pi_{\lambda,q}(x^{(3)} Z_{i} S^{\mathrm{op}}(x^{(2)}))\pi_{\lambda,q}(x^{(1)})) \otimes (\sigma_q(x^{(4)}) s_q( \psi_{j}) \sigma_q(S(x^{(5)}))) \sigma_q(x^{(6)})| \omega \rangle \\
&=& (\sum_{i,j} \alpha_{ij} \pi_{\lambda,q} (x'' \obad Z_{i}) \otimes s_q(x''' \ad \psi_{j}))(\pi_{\lambda,q}(x') \otimes \sigma_q(x''''))| \omega \rangle \\
&=& \epsilon(x'')D(\pi_{\lambda,q}(x') \otimes \sigma_q(x''')) | \omega \rangle = Dx.| \omega \rangle. \5 \square 
\end{eqnarray*}

The limit $q \rightarrow 1$ of this operator is the classical Dirac operator with a reduced connection. We can also add the cubical in $\mathrm{cl}_q(\g)$ to $D$ so that it satisfies the required properties.  The tensor product $V \otimes V$ contains an invariant subspace isomorphic to the adjoint representation in the negative spectral subspace of the braid operator. Thus, it does not vanish when embedded into the algebra  $\mathrm{cl}_q(\g)$. Let $\theta : V \rightarrow \mathrm{cl}_q(\g)$ denote the corresponding module isomorphism. The cubical part of the Dirac operator is defined by
\begin{eqnarray*}
\Gamma = 1 \otimes (m(\theta \otimes \gamma_q \circ \tau) (\Omega)) \in U_q(\g) \otimes \mathrm{cl}_q(\g)
\end{eqnarray*}
where $m$ is the product of $\mathrm{cl}_q(\g)$. $\Gamma$ is invariant under the action of $U_q(\g)$ and commutes with the representation. For a suitable choice of constant $N$ the operator $D + N \Gamma$ reduces to the geometric Dirac operator equipped with a Levi-Civita connection or to Kostant's Dirac operator \cite{Kos} in the limit $q \rightarrow 1$. 
\\

\textbf{10. Example: $SU_q(2)$.} For each $l \in \frac{1}{2} \mathbb{N}_0$ choose a basis of $2n+1$ dimensional vector space $V_l$ by $\{ |l,m \rangle: -l \leq m \leq l\}$. The irreducible finite dimensional representations $(V_{l}, \pi_{l,q}$) of $U_q(\mathfrak{su_2})$ are
\begin{eqnarray*}
\pi_{l,q}(k)|l,m \rangle &=& q^m |l,m \rangle  \\ \nonumber
\pi_{l,q}(e)|l,m \rangle &=& \sqrt{[l-m][l+m+1]} |l,m+1 \rangle \\ \nonumber 
\pi_{l,q}(f)|l,m \rangle &=& \sqrt{[l-m+1][l+m]} |l,m-1 \rangle.
\end{eqnarray*}
$U = V_{1/2}$ is the defining representation and $V = V_1$ is the adjoint representation.

The module isomorphism $\phi: V \rightarrow \mathfrak{L}^{\mathrm{op}}_q(\mathfrak{su_2})$ is defined by 
\begin{eqnarray*}
\phi(|1,1 \rangle) = t^{-1}e, \5 \phi(|1,0 \rangle) = \frac{1}{\sqrt{[2]}}(q^{-1}fe - qef),\5 \phi(|1,-1 \rangle) = - t^{-1}f.
\end{eqnarray*}
Let us write $\phi(|1,m \rangle) = Z_m$. 

Using the Clebsch-Gordan rule we find 
\begin{eqnarray*}
V \otimes V \simeq V_0 \oplus V \oplus V_2
\end{eqnarray*}
The subspace $V_0 \oplus V_2$ is in the positive spectral subspace of $\check{R}$ and the adjoint module $V$ is in the negative spectral subspace. The following Clifford algebra relations can be written down immediately by identifying the basis vectors of  $V_0 \oplus V_2$ with their image in $B_q$  
\begin{eqnarray*}
& &\psi_1 \psi_1 = \psi_{-1} \psi_{-1} = 0 \\
& & q^{-1} \psi_1 \psi_0 + q \psi_0 \psi_1 = 0 \\
& &q^{-2} \psi_1 \psi_{-1} + [2] \psi_0 \psi_0 + q^2 \psi_{-1} \psi_1 = 0 \\
& &\psi_0 \psi_{-1} + q^2 \psi_{-1} \psi_0 = 0 \\
& &\psi_1 \psi_{-1} + \psi_{-1} \psi_1 = b,
\end{eqnarray*}
where $\psi_i = \gamma_q(|1,i \rangle)$ and $b$ is some constant fixed from the normalization of the form $B_q$. 

The irreducible representation space for $\mathrm{cl}_q(\g)$ is $2$-dimensional, $\Sigma = V_{1/2}$. Let us choose $\sigma_q = \pi_{1/2,q}$. The vector space $B(\Sigma)$ becomes a $U_q(\g)$-module and the submodule isomorphic to the adjoint representation under the action \eqref{comp} has a basis
\begin{eqnarray*}
& &s_q(\psi_1) = \begin{pmatrix} 
  0 & \sqrt{q} \\
  0 & 0 \end{pmatrix}, \5 s_q(\psi_0) = - \frac{1}{\sqrt{[2]}}\begin{pmatrix} 
  q^{-1} & 0 \\
  0 & -q \end{pmatrix}, \5 s_q(\psi_{-1}) = \begin{pmatrix} 
  0 & 0 \\
  -\sqrt{q^{-1}} & 0 \end{pmatrix}
\end{eqnarray*}
These matrices satisfy the $q$-Clifford algebra relations with $b = -1$. 

The module isomorphism $\tau: V^* \rightarrow V$ is given by
\begin{eqnarray*}
\tau(|1,1 \rangle^*) = -q  [2]|1,-1 \rangle,\5 \tau(|1,0 \rangle^*) =  [2]|1,0 \rangle ,\5 \tau(|1,-1 \rangle^*) =  -q^{-1} [2]|1, 1 \rangle.
\end{eqnarray*}
Now we can write
\begin{eqnarray*}
D = \sum_l \phi(|1,l \rangle ) \otimes  \gamma_q \circ \tau( |1, l \rangle^*) = -q[2] Z_1 \otimes \psi_{-1} +[2] Z_0 \otimes \psi_0 - q^{-1}[2] Z_{-1} \otimes \psi_{1}.
\end{eqnarray*}
For any $l \in \frac{1}{2}\mathbb{N}_0$ we define $H_l = V_{l} \otimes \Sigma$. As a $U_q(\mathfrak{su(2)})$-module this decomposes as
\begin{eqnarray*}
H_l \simeq V_{l+1/2} \oplus V_{l-1/2},
\end{eqnarray*}
if $l > 0$ and $H_0 \simeq V_{1/2}$. Let us write $j\pm1/2 = j^{\pm}$. Denote by $|j^+,\mu \rangle$  $(-j^+ \leq \mu \leq j^+)$ and $|j^-,\mu \rangle$ $(-j^- \leq \mu \leq j^-)$ the basis vectors of the irreducible components. The Dirac operator acts on $H_l$ by 
\begin{eqnarray*}
D.|j^+ \mu \rangle = [2j] |j^+ \mu \rangle,\5 D.|j^- \mu \rangle = -[2j+2] |j^- \mu  \rangle 
\end{eqnarray*}
if $l > 0$ and $H_0 \in \mathrm{Ker}(D)$. \\

\textbf{11.} The matrix elements of the irreducible finite dimensional representations of $U_q(\mathfrak{su_2})$ span the space of polynomial functions on the quantum group 
\begin{eqnarray*}
\mathbb{C}[SU_q(2)] = \bigoplus_{l \in \frac{1}{2}\mathbb{N}_0} V_l^* \otimes V_l \simeq  \bigoplus_{l \in \frac{1}{2}\mathbb{N}_0} V_l \otimes V_l.
\end{eqnarray*}
The multiplication is derived from the Clebsch-Gordan coefficients. We use the Haar state of $\mathbb{C}[SU_q(2)]$ to complete $\mathbb{C}[SU_q(2)] \otimes \Sigma$ to a Hilbert space $\h$. The algebra $U_q(\g)$ acts from left by 
\begin{eqnarray*}
x.(|l,m \rangle \otimes |l,n \rangle \otimes |\frac{1}{2},s \rangle) = (\pi_{l,q}(x')|l,m \rangle) \otimes |l,n \rangle \otimes  (\pi_{1/2,q}(x'') |\frac{1}{2},s \rangle), 
\end{eqnarray*}
for all $|l,m \rangle \otimes |l,n \rangle \otimes |\frac{1}{2}, s \rangle \in \h$. The explicite decomposition of the prehilbert space into irreducible components under the left action is given in \cite{DB} 
\begin{eqnarray*}
(\bigoplus_{l \in \frac{1}{2} \mathbb{N}_0} V_{l} \otimes V_{l}) \otimes \Sigma \simeq V_{1/2} \oplus \bigoplus_{l \in \frac{1}{2} \mathbb{N}} (V_{j+1/2} \otimes V_{j}) \oplus (V_{j-1/2} \otimes V_{j}) := W^{\uparrow}_0 \oplus \bigoplus_{l \in \frac{1}{2} \mathbb{N}} W^{\uparrow}_j \oplus W^{\downarrow}_j.
\end{eqnarray*}
The components $W^{\uparrow}_j$ and $W^{\downarrow}_j$ have multiplicities $(2j+2)(2j+1)$ and $2j(2j+1)$. The orthonormal basis of $\h$ is chosen by
\begin{eqnarray*} 
 |j \mu n \uparrow \rangle \in W^{\uparrow}_j, \5 |j' \mu' n \downarrow  \rangle \in W^{\downarrow}_j 
\end{eqnarray*}
for $j \in \frac{1}{2} \mathbb{N}_0$, $ j' \in  \frac{1}{2} \mathbb{N}$, $ -j^+ \leq \mu \leq j^+$, $ -j^- \leq \mu' \leq j^-$ and $ -j \leq n \leq j$.

Let us adopt the column vector notation of \cite{DB}
\begin{eqnarray*}
|j \mu n \rangle \rangle = \begin{pmatrix} 
  |j \mu n \uparrow \rangle \\
  |j \mu n \downarrow \rangle
  \end{pmatrix}.
\end{eqnarray*}
A faithful $*$-representation for the algebra $\mathbb{C}[SU_q(2)]$ on $\h$ is developed in \cite{DB} (Proposition 4.4), where the notation matches with ours. This representation coincides with the GNS representation. Since we know how $D$ acts on each irreducible piece we get
\begin{eqnarray*}
D.|j \mu n \rangle \rangle = \begin{pmatrix} 
  [2j] & 0 \\
  0 & -[2j+2] \end{pmatrix} |j \mu n \rangle \rangle,
\end{eqnarray*}
when $j > 0$ and $W^{\uparrow}_0 \in \mathrm{Ker}(D)$. The kernel is nontrivial and we define an approximated sign operator in the usual way: $F = D(1+D^2)^{-1/2}$. \\

\noindent \textbf{Proposition.} Triple $(\mathbb{C}[SU_q(2)], F, \h)$ is a $1$-summable Fredholm module. \\

\noindent Proof. The spectrum of $D$ grows exponentially as a function of $j$ and therefore $F^2-1$ is a trace class operator
\begin{eqnarray*}
F^2-1 = -(1+D^2)^{-1} \in L^1(\h). 
\end{eqnarray*}
The operators of the irreducible $*$-representation \cite{DB} of $\mathbb{C}[SU_q(2)]$ on $\h$ are  of the form
\begin{eqnarray*}
\eta(x) = \sum_{i} X(i),\5 X(k) = \begin{pmatrix}  
  X(k)_{\uparrow \uparrow} & X(k)_{\uparrow \downarrow} \\
  X(k)_{\downarrow \uparrow} & X(k)_{\downarrow \downarrow} \end{pmatrix}, 
\end{eqnarray*} 
where $X(k)_{\uparrow \uparrow}, X(k)_{\downarrow \downarrow} \in B(\h)$ and $X(k)_{\uparrow \downarrow}, X(k)_{\downarrow \uparrow} \in L^1(\h)$ and  each $X(k)$ shifts the index $j$ of the basis vector by $k \in \frac{1}{2} \mathbb{Z}$. The sum is always finite. We have 
\begin{eqnarray*}
[F, X(k)_{\uparrow \uparrow}]|j \mu n \rangle = \Big( \frac{[j+k]}{(1 + [j+k]^2)^{1/2}} - \frac{[j]}{(1 + [j]^2)^{1/2}} \Big) X(k)_{\uparrow \uparrow}|j \mu n \rangle.
\end{eqnarray*}
The first term can be easily seen to give a sequence which decays rapidly as a function of $j$. Therefore $[F, X(k)_{\uparrow \uparrow}] \in L^1(\h)$. Similarly we see that $[F, x_{\downarrow \downarrow}] \in L^1(\h)$. The commutators $[F, X(k)_{\uparrow \downarrow}]$ and $[F, X(k)_{\downarrow \uparrow}]$ are in trace class because the off diagonal blocks are trace class operators and $F$ is bounded. Thus, $[F,X(k)] \in L^1(\h)$ for all $k$ and thus $[F, \eta(x)] \in L^1(\h)$. $\5 \square$ \\

\textbf{12. Acknowledgements.} The author wishes to thank his supervisor Jouko Mickelsson for several helpful discussions. This project was supported by the V$\ddot{a}$is$\ddot{a}$l$\ddot{a}$ foundation of the Finnish Academy of Science and Letters and The Finnish National Graduate School in Mathematics and its Applications.

\end{document}